\documentclass[11pt,reqno]{amsart}
\usepackage[latin1]{inputenc}
\usepackage[english]{babel}
\usepackage{mathrsfs}
\usepackage{array}
\usepackage{multicol}
\usepackage[dvips]{epsfig}
\usepackage{amsmath}
\usepackage{graphics}
\usepackage{color}
\usepackage{latexsym}
\usepackage{amsmath}
\usepackage{amssymb}
\usepackage[dvips]{epsfig}
\newtheorem{Th}{Theorem}[section]
\newtheorem{Le}{Lemma}[section]
\newtheorem{co}{Corollary}[section]
\newtheorem{de}{Definition}[section]

\newtheorem{re}{Remark}[section]

\newcommand{\dsps}{\displaystyle}

\newcommand{\1}{1\hspace{-1.4mm}1}

\begin{document}

\title[Nonlocal Second-order equations]{Nonlocal Second-order geometric equations arising in tomographic reconstruction}
\author[ali srour]{ ali srour} 
\address{Laboratoire de Mathématiques et Physique Théorique
Université François-Rabelais Tours, Fédération Denis Poisson -UMR CNRS 6083,
Parc de Grandmont, 37200 Tours. France.
 Laboratoire de Mathématiques et Applications  
Physique Mathématique d'Orléans, 45067 Orléans cedex 2} 

\date{\today}
\email{srour@lmpt.univ-tours.fr }
\keywords{Nonlocal Hamilton-Jacobi equations, Tomographic reconstruction, Nonlocal front propagation, level-set approach, Viscosity solutions}
\begin{abstract}
In this paper, we study a new model of nonlocal geometric   equations which appears in tomographic reconstruction when using the level-set method. We treat two additional difficulties which make  the work original. On  one hand, the level lines do not evolve along  normal directions, and  the nonlocal term is not of  ``convolution type". On the other hand,  the speed is not necessarily bounded compared to the nonlocal term. We  prove a existence and uniqueness results of our model.
\end{abstract}
\maketitle 
In this paper, we study a fully nonlinear parabolic equation with nonlocal term. More precisely, 
\begin{equation}\label{HG}
\begin{cases}
\displaystyle\frac{\partial u}{\partial t}(x,t)= F(x,t,Du,D^{2}u,K^{k,+}_{{x,t,u}} )\hspace{0.5cm} \rm{in}~~\mathbb{R}^  {N} \times \left[0,T\right],\\
 \\
 u(x,0)=u_{0}(x) \hspace{0.5cm}{\rm{in}}~~ \mathbb{R}^{N}\times  \left\{0\right\},
 
\end{cases}
\end{equation}
where $N\geq 1$ is an integer, $T>0$. The unknown function is   $u:\mathbb{R}^{N}\times[0,T]\longrightarrow \mathbb{R}$, $Du$ and $D^{2}u$ denote respectively the gradient and the Hessian of $u$ with respect to the space variable, $u_{0}$ : $\mathbb{R}^{N}\longrightarrow \mathbb{R}$ is the given initial data and $K^{k,+}_{x,t,u}$ denotes the nonlocal term,  given for  some integer $k \leq  N$ by
\begin{center}
$$K^{k,+}_{{x,t,u}}=[u]^{+}_{x,t}\cap A^{k}_{x}$$
\end{center}
where $$[u]^{+}_{x,t}=\left\{y\in \mathbb{R}^{N} : u(y,t)\geq u(x,t)\right\}$$ and $$A^{k}_{x}=\left\{y=(y_{1},...,y_{k},y_{k+1},...,y_{N})\in \mathbb{R}^{N}:\pi_{k}(y)=\pi_{k}(x)\right\}$$
where $\pi_{k}$ denote  the  projection  function from $\mathbb{R}^{N}$ into $\mathbb{R}^{k}$ defined by $$\pi_{k}(x)=(x_{1},...,x_{k})$$ 
for all $x=(x_{1},...,x_{k},x_{k+1},...,x_{N})$. In the same manner 
$$K^{k,-}_{{x,t,u}}=[u]^{-}_{x,t}\cap A^{k}_{x}$$
where $$[u]^{-}_{x,t}=\left\{y\in \mathbb{R}^{N} : u(y,t)> u(x,t)\right\}.$$ 
The  nonlinearity  $F$ is a continuous function from $\mathbb{R}^{N}\times\mathbb{R}\times \mathbb{R}^{N}\backslash\{0\}\times\mathscr{S}_{N}\times \mathscr{B}_{N-k}$ into $\mathbb{R}$, where $\mathscr{S}_{N}$ is the set of real symmetric $N\times N$ matrices and  $\mathscr{B}_{N-k}$ is the set of  equivalence classes of all subsets of $\mathbb{R}^{N-k} $ with respect the relation $A\sim B$ if $\mathcal{L}^{N-k}(A\Delta B)=0$, where $\mathcal{L}^{N-k}$ is the Lebesgue measure on $\mathbb{R}^{N-k}$. 

We consider $\mathscr{B}_{N-k}$ with a  topology that comes from the metric $$\dsps d(A,B)=\sum^{\infty}_{n=0}\frac{\mathcal{L}^{N-k}\Big((A\Delta B) \cap B(0,n)\Big)}{2^{n}\mathcal{L}^{N-k}(B(0,n))},$$ where $B(0,n)$ denotes  the ball in $\mathbb{R}^{N-k}$ of center $0$ and radius $n$. With this topology a sequence $(K^{k}_{n})_{n\geq1}$ in  $\mathscr{B}_{N-k}$ converges to $K^{k}\in\mathscr{B}_{N-k}$ if and only if $\dsps\1_{K^{k}_{n}}$ converges to $\1_{K^{k}}$ in $L^{1}_{loc}(\mathbb{R}^{N-k})$. 
 
We are interested in equations of type (\ref{HG}) which are related with tomographic reconstruction using active curves and the level-set approach \cite{IR02,K}. The model case we have in mind is 
 \begin{equation}\label{Te}
F(x,t,p,X,K^{1,+}_{{x,t,u}})=C(x,t)\left(\int_{K^{1,+}_{{x,t,u}}} g(z)dz\right)|p|+{\rm{Trace}}\left[(I-\frac{p\otimes p}{|p|^{2}})X\right],
\end{equation}
where $x=(x_{1},x_{2})$, $K^{1,+}_{{x,t,u}}= [u]^{+}_{x,t}\cap A^{1}_{x}$, $g$ is a positive function  and $C$ is  positive and Lipschitz continuous function. 
 
We recall that the level-set approach was first introduced by Osher and Sethian \cite{osher88} for numerical computations, and then developed from  a theoretical point view by Evans and Spruck \cite{ES91} for motion by mean curvature and by Chen, Giga and Goto \cite{GGC91} for general normal velocities. We also refer the reader to Barles, Soner and Souganidis \cite{BSS93} and Souganidis  \cite{S95,S97} for different presentations and other results on the level-set approach.

In \cite{SL03} Slep\v{c}ev studied the motion of fronts in bounded domains by normal velocities  which  can depend on non local terms, in addition to the curvature, the normal direction and the location of the front. In fact, the velocities depend on non local terms, if the velocities at any point of the front depend on the set that the front encloses. Depending on  the velocities, the motion of the front can be described by the partial differential equation
\begin{equation}\label{SL}
\dsps\frac{\partial u}{\partial t}(x,t)=F(x,t,Du(x,t),D^{2}u(x,t),[u]^{+}_{x,t}) \;\;\; {\rm{in}}\quad O \times [0,T]
\end{equation}
\noindent with the Neumann boundary conditions $\dsps\partial u/\partial\gamma=0$~~${\rm{on}}$ ~$\partial O \times[0,T] $ and $u(x,0)=u_{0}(x)$, where $O$ is a bounded domain in $\mathbb{R^{N}}$. Slep\v{c}ev proved an existence and uniqueness result for this equation, using the viscosity solution.

In \cite{BCLM}, Barles,  Cardaliaguet,  Ley,  Monneau studied the first-order nonlocal equation
\begin{equation}
\begin{cases}
\displaystyle\frac{\partial u}{\partial t}(x,t)=\left(C_{0}(.,t)*\1_{[u]^{+}_{x,t}}(x,t)+C_{1}(x,t)\right)|Du| \hspace{0.5cm} {\rm{in}}~~\mathbb{R}^{N} \times \left[0,T\right],\\
 \\
 \displaystyle u(x,0)=u_{0}(x) \hspace{0.5cm} {\rm{in}}~~\mathbb{R}^{N}\times  \left\{0\right\}. 
\end{cases}
\end{equation}
This equation appears when modelling dislocations in crystals using the level set approach. The $\ast$ denotes the convolution in space,  $C_{0}$ and  $C_{1}$ are two functions on which  we have some conditions. They proved an existence and uniqueness result for this equation in  both cases where  $C_{0}$ is a positive and a negative  function, $C_{0}(.,t)\in L^{1}(\mathbb{R}^{N})$ for all $t\in [0,T]$  and $C_{0}, C_{1}$  are two Lipschitz continuous functions .

Now a  question  comes out: what are   the differences between our equation $(\ref{HG})$ and the equations $(3)$ and $(4)$? To explain  the differences between our equation (\ref{HG}) and the equation (\ref{SL}) for example, we consider the case $ N=2$ and the typical equation (\ref{Te}).

In $(4)$, for any $(x,t)$$\in \mathbb{R}^{2}\times [0,T]$, the nonlocal term is given by
\begin{center}
$[u]^{+}_{x,t}=\left\{y\in \mathbb{R}^{2}:u(y,t)\geq u(x,t)\right\}.$
\end{center}
and we integrate  on  a subset of $\mathbb{R}^{2}$. In our Equation $(2)$, the  new nonlocal term is given  for any $(x,t)$$\in \mathbb{R}^{2}\times [0,T] $ by 
\begin{center}
$K^{1,+}_{{x,t}}=[u]^{+}_{x,t}\cap A^{1}_{x}= \left\{y\in \mathbb{R}^{2}:u(y,t)\geq u(x,t)\right\}\cap\left\{y_{1}=x_{1}\right\}$
\end{center}
and we integrate in $(2)$ over a straight line i.e. the second variable is always fixed. We remark that the technique used in \cite{BCLM,SL03}  cannot be applied in the case of  this new nonlocal term. We  change the dependence in $\dsps Du $, we  are able to prove a uniqueness and existence result of our equation $(\ref{HG})$.

Moreover, contrary to the cases studied in \cite{BCLM,SL03}, in the initial compact front, we allows a not bounded dependence of  volume in our equation  (\ref{HG})(see \textbf{(H5-1)} and \textbf{(H6-1)}).

Let us now explain how this paper is organized: in Section $2$, we present our model and  we recall the definition of viscosity solutions. In Section $3$, we prove a uniqueness result for our model with compact fronts and  in Section $4$ we show a uniqueness result for non compact fronts. In Section $5$, we  prove an existence result for compact and non compact fronts. Finally, we give, in Appendix, the proof of the stability result for our model.
\section{ Presentation of our model  }
 Tomography processes are widely studied by many authors. Most of cases concern the images restoration from a lot of projected data. In \cite{MS}, M. Somekh addresses the problem of reconstruction of an image from two pairs of its orthogonal projections. The paper by Dinten, Bruandet and Peyrin  \cite{K} addresses tomographic reconstruction of  binary objects from a small number of noisy projections in applications where the global dose remains constant with an increase or a decrease of the number of projections. Here we are  interested in   the reconstruction  method  of \cite{IR02}, where the authors consider a  single view of tomographic reconstruction for radially symmetric objects and binary image.  They formulate  the  problem as a  front propagation  which consists in evolving the contour of the noisy image (initial front),  with a selected normal velocity, so that it converges toward  contour of the initial object. This evolution is described by an approach of level-set which leads a nonlocal equation      

\begin{equation}\label{F}
\begin{cases}
\dsps\frac{\partial u}{\partial t}=C(x,t)\left(\int_{\mathbb{R}} \1_{\{u(.,x_{2},t)\geq 0\}}(z) g(z)dz\right)|Du|+ {\rm{Trace}}\left[(I-\frac{Du\otimes Du}{|Du|^{2}})D^{2}u\right],
 \\
 u(x,0)=u_{0}(x) \hspace{0.5cm}{\rm{in}}~~ \mathbb{R}^{N}\times  \left\{0\right\},
 
\end{cases}
\end{equation}
where $x\in \mathbb{R}^{2}$, $g $ is  a positive functions and $g \in {L}^{1}(\mathbb{R})$. The second term  on the     right hand  side of (\ref{F}) is the mean curvature. This term was studied  in \cite{ES91} where it  was used to regularize the evolution of the initial front with normal velocity. This equation (\ref{F}) has been studied in \cite {IR02} with $C\leq0$, but in this papers we change the sign of $C$ in order to keep the assumption  of monotony (see (\textbf{\textbf{H3}})) which is a classical assumption to provide  comparison and uniqueness  results. From this equation, since $x_{2}$ is fixed, a discontinuity appears in the normal speed of  propagation. To see this problem of continuity  it is enough to start with a rectangular initial front whose sides are parallel to the coordinate axes and this rectangle can be represented as a level-set zero of a function $u$ from $\mathbb{R}^{N}\times[0,T]$ to $\mathbb{R}$ i.e $u=0$ on the sides, $u>0$ inside the rectangle and $u<0$ outside (see figure 1). We suppose that this front evolve with a normal simplified velocity given by 
$$v(x_{1},x_{2},t)=\int_{K=\{s\in\mathbb{R}: u(s,x_{2},t)\geq0\}}ds.$$

\begin{figure}[!ht]
\begin{center}
\epsfig{file =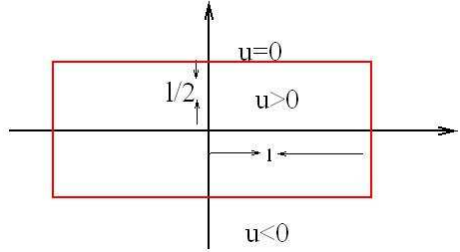,scale=0.5}
		\caption{rectangular initial front}	
\end{center}
\end {figure}
We consider now the points $\dsps(x,\frac{l}{2},t)$ where $l$ is the width of the rectangle and $2l$ his longer. It is easy to prove that  the velocity  on these points is $2l$, but if we move vertically to $\dsps(x,\frac{l}{2}+\epsilon)$, the velocity becomes $\dsps\dsps v(x,\frac{l}{2}+\epsilon,t)=0.$\\
Then, if we move vertically, a jump appears in the velocity. In the present work, instead of evolving in the normal direction, we move in the horizontal direction. Thus with the horizontally velocity, the horizontals sides of the rectangle remains fixed, whereas the vertically sides move in the horizontal direction but in two different sens. After this modification, the equation (\ref{F}) becomes
\begin{equation}\label{F'}
\begin{cases}
\dsps\frac{\partial u}{\partial t}=C(x,t)\left(\int_{\mathbb{R}} \1_{\{u(.,x_{2},t)\geq 0\}}(z) g(z)dz\right)|\frac{\partial u}{\partial x_{1}}|+ {\rm{Trace}}\left[(I-\frac{Du\otimes Du}{|Du|^{2}})D^{2}u\right],
 \\
 u(x,0)=u_{0}(x) \hspace{0.5cm}{\rm{in}}~~ \mathbb{R}^{N}\times  \left\{0\right\}.
 \end{cases}
 \end{equation}
As it was remarked by Slep\v{c}ev \cite{SL03}, in the level-set approach, all level-sets of the solution $u$ should have the same type of normal velocity. A nonlocal term using $[u]^{+}_{x,t}$ instead  of $\{u(.,t)\geq 0\}$ is more appropriate. It is why we consider  
\begin{equation}\label{SF}
\dsps\frac{\partial u}{\partial t}=C(x,t)\left(\int_{\mathbb{R}} \1_{\{u(.,x_{2},t)\geq u(x_{1},x_{2},t)\}} g(z)dz\right)|\frac{\partial u}{\partial x_{1}}|+ {\rm{Trace}}\left[(I-\frac{Du\otimes Du}{|Du|^{2}})D^{2}u\right].
\end{equation}
With this modification,  we will prove in Section 3  uniqueness and existence results  for (\ref{SF}).

We now list the basic requirements on $F$. We point out that the main assumption introduced because of the presence of the nonlocal term is the monotonicity with respect to set inclusion and the continuity of $F$ with respect to the topology  on $\mathscr{B}_{N-k}$ which is defined in the introduction.

\vspace{3mm}

\noindent\textbf{(H1)} $F$ is continuous from $\mathbb{R}^{N}\times\mathbb[0,T]\times \mathbb{R}^{N}\backslash\{0\}\times\mathscr{S}_{N}\times \mathscr{B}_{N-k}$ into $\mathbb{R}.$

\vspace{3mm}

\noindent\textbf{(H2)} For all $(x,t,K^{k},L^{k})\in \mathbb {R}^{N} \times [0,T]\times \mathscr{B}_{N-k}\times\mathscr{B}_{N-k}$, we have $$- \infty <F_{\star}(x,t,0,O,K^{k})=F^{\star}(x,t,0,O,L^{k})<+\infty .$$
Here $F^{\star}$ denotes the upper semicontinuous envelope of $F$, while $F_{\star}$ is the lower semicontinuous envelope of $F$. These functions are defined by
 $$ F^{\star}(x,t,p,X,K^{k})= \dsps\mathop{{\rm{lim sup}}}_{n}{F(x_n,t_n,p_n,X_n,K^{k}_n)},$$
$$ F_{\star}(x,t,p,X,K^{k})= \dsps\mathop{{\rm{lim inf}}}_{n}{F(x_n,t_n,p_n,X_n,K^{k}_n)},$$
where $(x_n,t_n,p_n,X_n)$ converges to $(x,t,p,X)$ in $\mathbb{R}^{N}\times[0,T]\times\mathbb{R}^{N}\backslash\{0\}\times\mathscr{S}_{N}$ and $K^{k}_n$ convergs to   $K^{k}$ with respect the topology in $ \mathscr{B}_{N-k}.$

\vspace{3mm}

\noindent\textbf{(H3) }$F$ is nondecreasing in its set argument i.e, for any $K^{k}$, $L^{k}$ in $\mathscr{B}_{N-k}$ such that $K^{k}\subset L^{k}$, we have 
$$ F(x,t,p,X,K^{k})\leq F(x,t,p,X,L^{k})$$   
for any $X$ in $\mathscr{S}_{N}$, and for any $(x,t,p)$ in $\mathbb{R}^{N}\times[0,T]\times \mathbb{R}^{N}\backslash\{0\}.$

\vspace{3mm}

\noindent\textbf{(H4)}  $F$ is geometric: for any $\lambda >0$, $\mu \in \mathbb{R}$, we have
 $$F(x,t,\lambda p,\lambda X+\mu p\otimes p,K^{k})=\lambda F(x,t,p,X,K^{k})$$
 for all $(x,t,p,X,K^{k}) \in \mathbb{R}^{N}\times[0,T]\times \mathbb{R}^{N}\backslash\{0\}\times\mathscr{S}_{N}\times\mathscr{B}_{N-k}.$
\begin{de}\label{ds} {\rm{ \textbf{(Slep\v{c}ev \cite{SL03}}})} \\ An upper-semicontinuous function $u$ : $\mathbb{R}^{N}\times[0,T]\longrightarrow \mathbb{R}$,   is a viscosity subsolution of $(\ref{HG})$  if, for any $\phi \in C^{2}(\mathbb{R}^{N}\times[0,T])$, for any  maximum point $(\bar{x},\bar{t})$ of $u-\phi$, if  $\bar{t} > 0$  then
$$\frac{\partial \phi}{\partial t}(\bar{x},\bar{t})\leq F_{\star}(\bar{x},\bar{t},D\phi(\bar{x},\bar{t}),D^{2}\phi(\bar{x},\bar{t}) ,K^{k,+}_{\bar{x},\bar{t},u})$$
and  $u(\bar{x},0)\leq u_{0}(\bar{x})$ if $\bar{t}=0.\\$
A lower-semicontinuous function $u$ : $\mathbb{R}^{N}\times[0,T]\longrightarrow \mathbb{R}$ is a viscosity supersolution of $(\ref{HG})$ if, for any $\phi \in C^{2}(\mathbb{R}^{N}\times[0,T])$, for any  minimum point $(\bar{x},\bar{t})$ of $u-\phi$, if $\bar{t} > 0 $ then
$$\frac{\partial \phi}{\partial t}(\bar{x},\bar{t})\geq F^{\star}(\bar{x},\bar{t},D\phi(\bar{x},\bar{t}),D^{2}\phi(\bar{x},\bar{t}) ,K^{k,-}_{\bar{x},\bar{t},u})$$
and $u(\bar{x},0)\geq u_{0}(\bar{x})$ if $\bar{t}=0$.\\
A function is a viscosity solution of $(\ref{HG})$, if it is both a subsolution and supersolution of  $(\ref{HG})$.
\end{de}
\begin{re}\label{rs}
\rm{In the definition of a subsolution and a supersolution, ``test sets'' were chosen differently. This is a major point when viscosity solutions are to be extended to non local, geometric parabolic equations. If  ``$\geq$'' were used instead of ``$>$'' in the definition of supersolutions, existence results, among other things, would hold no more. For more  information, we refer the reader to the arguments given by Slep\v{c}ev in \cite[Definition 2.1]{SL03}. For convenience,  we proved another  proof for these arguments in  the Appendix.
}\end{re}
In the initial compact front  we will seek the solutions of our equation (the solution of which represents the front via its zero level-set) by functions in the  class $\mathcal{C}$  given by the following definition.
 
\begin{de}  A  function $u$ : $\mathbb{R}^{N}\times[0,T]\rightarrow \mathbb{R}$ is in the class $\mathcal{C}$, if $-1<u(x,t)\leq 1$ and  $\dsps\mathop{\rm{inf}}_{x\in \mathbb{R}^{N}}{u(x,t)}=\dsps\mathop{\lim}_{|x|\rightarrow+\infty}{u(x,t)}=-1$ for all $x$ {\rm{in}} $ \mathbb{R}^{N}$ and  uniformly with respect to $t\in[0,T]$.
\end{de} 
\begin{re}\rm{ All continuous functions in $\mathcal{C}$ are  uniformly continuous bound-ed function and  any compact front can  be  represented as a level-set zero of  a function in $\mathcal{C}$. Indeed, let $E\subset \mathbb{R}^{N}$ be a bounded open set and $\Gamma_{0}=\partial E$ be the initial compact front. This front is represented by the signed function defined by 
$$ d^{s}_{\partial E}= d(x,\mathbb{R}^{N}\backslash{E})-d(x,E),$$ 
where $d(x,E)$ denote the distance function to $E$. It is easy to check that $ d^{s}_{\partial E}> 0$ if $x\in E$, 
$d^{s}_{\partial E}< 0$ if $x \in \mathbb{R}^{N}\backslash\bar{E}$ and  $d^{s}_{\partial E}= 0$ if and only  if $x\in \Gamma_{0}$. Moreover $d^{s}_{\partial E}$ is Lipschitz continuous with constant 1. Now, we consider the function $\dsps f(x)=2\;{\rm{arctan}}(d^{s}_{\partial E}(x))/\pi$ and we use the properties of the arctan function to prove that $f(x)>-1$ for all $x\in \mathbb{R}^{N}$ and $f(x)$ tends to $-1$ as $x$ tends to $\infty$. Then, $f\in \mathcal{C}$ and the initial compact front $\Gamma_{0}$
is represented by $f$. A typical example of  function which belongs to this  class is  $\dsps g(x,t)=\frac{e^{-At}}{1+|x|^{2}}-1$ where $(x,t)\in\mathbb{R}^{N}\times[0,T]$ and  $A$ is a positive constant.}\end{re}
Now, another question comes out: for which reason, we use the class $\mathcal{C}$  for sub and supersolutions in the compact front case? This class  $\mathcal{C}$ allows  to treat the sets $K^{k,\pm}_{x,t,u}$ of infinite volume. 
 To explain that, we start with an initial front which propagates with the following velocity 
$$ V(x,t,[u]^{+}_{x,t})=\dsps\left(\int_{[u]^{+}_{x,t}}dz\right)^{\frac{1}{N}}=\left({\rm{vol}}\{u(.,t)\geq u(x,t)\}\right)^{\frac{1}{N}}.$$
When we describe this evolution by the level-set method we obtain the nonlocal level-set equation 
\begin{equation}\label{EV}
\begin{cases}
\dsps\frac{\partial u}{\partial t}(x,t)=\left({\rm{vol}}\{u(.,t)\geq u(x,t)\}\right)^{\frac{1}{N}}|Du|\hspace{0.5cm} {\rm{in}}~~\mathbb{R}^  {N} \times \left[0,T\right],\\
\displaystyle u(x,0)=u_{0}(x) \hspace{0.5cm}{\rm{in}}~~ \mathbb{R}^{N}\times  \left\{0\right\} .
\end{cases}
\end{equation}
First, the power  $ \frac{1}{N}$  ensures  the existence of a viscosity solution of (\ref{EV}) and  ensures  that the front does not explode on short  time $t$ (see \cite[Remark 4.1]{BL}). Moreover, if the solution $u$ belongs to the class $\mathcal{C}$, for any $(x,t)\in\mathbb{R}^{N}\times[0,T]$ there exists $M_{x,t}$ such that, the Lebesgue measure $\mathcal{L}^{N}([u]_{x,t})\leq M_{x,t}< \infty$. Whence, the  previous definition of viscosity solution remains valid for the equation (\ref{EV}).
\section{Uniqueness result for compact fronts}
In the case of compact fronts, we start with a compact initial front and we  prove that the front remains compact i.e.  there exists a solution $u$ of (\ref{HG}) which lies in $\mathcal{C}$. We use the following assumptions:

\vspace*{3mm}

\noindent\textbf{(H0)} \textbf{Initial compact front}: $u_{0} \in {\rm{C}}(\mathbb{R}^{N})$ and  $u_{0}\in\mathcal{C}$.\\
Now, for $\theta \in\{0,1\}$, we have the following assumptions.

\vspace*{3mm}

\noindent\textbf{(H5-$\theta$)} There exist positive constants $L_{1}$ and $L_{2}$ such that for $(x,t,p,X,K^{k})$ in $\mathbb{R}^{N}\backslash\{0\}\times[0,T]\times \mathbb{R}^{N}\times \mathscr{S}_{N}\times\mathscr{B}_{N-k}$, we have
$$|F(x,t,p,X,K^{k})|\leq \dsps L_{1}\Big((1+|x|)|p|+\dsps\theta\mathcal{L}^{N-k}(K^{k})^{\frac{1}{N-k}}|\tilde{\pi}_{k}(p)|\Big)+\dsps L_{2}(1+|x|^{2})|X|_{\infty} $$
where $|A|_{\infty}= \dsps\mathop{\rm {max}}_{ 1 \leq i\leq N}{\dsps \sum^{N}_{j=1}|a_{i,j}}|$, for  any matrix $A=(a_{i,j})_{1\leq i,j\leq N} ~\rm{in} ~\mathscr{S}_{N}$ and $\tilde{\pi}_{k}$ denote the projection from $\mathbb{R}^{N}$ into $\mathbb{R}^{N-k}$ defined by $\tilde{\pi}_{k}(x)=(x_{k+1},...,x_{N})=x-\pi_{k}(x),$ for all $x\in \mathbb{R}^{N}$.

\vspace{3mm} 

\noindent\textbf{(H6-$\theta$)}  There exists a nondecreasing modulus of continuity $w :[0,\infty] \rightarrow [0,\infty]$ which  satisfies $w(0+)=0 $   and 
$$F(x,t,p,X,K^{k})-F(y,t,p,Y,K^{k})\leq \dsps w \bigg(|x-y|\Big(1+|p|+\theta\mathcal{L}(K^{k})^{\frac{1}{N-K}}|\tilde{\pi}_{k}(p)|\Big)+\sigma |x-y|^{2}\bigg)$$
whenever $(x,y,t,p,K^{k} )\in \mathbb{R}^{N}\times\mathbb{R}^{N}\times[0,T]\times\mathbb{R}^{N}\backslash\{0\}\times\mathscr{B}_{N-k}$, $\sigma$ is a positive constant and $X,Y$ satisfying  the inequality 
$
\left(
\begin{array}{cc}
X    &   0 \\

0    & -Y   \\
\end{array}
\right)
$
$\leq $
$
\left(
\begin{array}{cc}
Z   &   -Z \\

-Z   & Z   \\
\end{array}
\right),
$
 where $Z \in \mathscr{S}_{N}$  and $\leq$ stands for the partial ordering in $\mathscr{S}_{N}$.

\vspace{3mm}
 
\noindent\textbf{(H7)} For $1<k \leq N$, for all $p$ in $\mathbb{R}^{N}$, there exists a  continuous function $G$ from $\mathbb{R^{N}}\times[0,T]\times \mathbb{R}^{k}\backslash\{0\}\times \mathscr{S}_{N}$ to $\mathbb{R}^{N-k}$  such that: if $\tilde{\pi}_{k}(p)=0$ then  
$$F(x,t,p,X,K^{k})=G(x,t,\pi_{k}(p),X)$$
for all $(x,t,p,X) \in \mathbb{R}^{N}\times[0,T]\times \mathbb{R}^{N}\backslash\{0\}\times\mathscr{S}_{N}$. Moreover, we suppose that, $G$  satisfies the assumption $\textbf{(H6-0)}$.
\begin{Th}\label{th1}Assume {\rm{\textbf{(H0)-(H1)-(H2)-(H3)-(H4)-(H5)-(H6-1)}}} and \textbf{(H7)}. Let $u \in\mathcal{C}$  (resp. $v\in \mathcal{C}$) be a bounded upper-semicontinous subsolution of (\ref{HG}) (resp.  bounded lower-semicontinous supersolution of (\ref{HG})), then $u\leq$ $v$ in $\mathbb{R}^{N}\times [0,T]$.
\end{Th}
\begin{co}\label{CO} Under the assumptions of Theorem $\ref{th1}$, there exists a unique viscosity solution $u\in \mathcal{C}$ of $(\ref{HG})$.
\end{co}
\noindent The  proof of this Corollary is postponed to (Section 5).

\vspace{3mm}

\noindent\textbf{Proof of  Theorem \ref{th1}.}\\
\textbf{1. The test-function}. We argue by contradiction assuming there exists $(\hat{x},\hat{t})\in \mathbb{R}^{N}\times[0,T]$  such that $$u(\hat{x},\hat{t})-v(\hat{x},\hat{t})=M>0.$$
Since $u$ and $v$ are bounded, the following supremum 
\begin{equation}\label{max}
\displaystyle M_{\epsilon,\eta}= \mathop{\rm {sup}}_{(x,y,t)\in (\mathbb{R}^{N})^{2}\times [0,T]}\{{u(x,t)-v(y,t)-\displaystyle \sum^{N}_{i=0}\frac{|x_{i}-y_{i}|^{4}}{4\epsilon^{2}}-\eta t }\}
\end{equation}
is finite for any $ \epsilon,\eta  > 0 $. We choose  $\eta $ is  small enough so that 
\begin{equation}\label{maxp}
M_{\epsilon,\eta}\geq \frac{M}{2}>0.
\end{equation}
Since $u,v \in \mathcal{C}$, we have 
$$\mathop{\lim}_{|x|\rightarrow +\infty}[{u(x,t)-v(x,t)}]=0$$  uniformly with respect to $t\in[0,T]$. Therefore from (\ref{maxp}), the   supremum in $(\ref{max})$ is achieved at a point $(\bar{x},\bar{y},\bar{t}),$ 
\begin{equation}\label{S}
M_{\epsilon,\eta}=u(\bar{x},\bar{t})-v(\bar{y},\bar{t})-\sum^{N}_{i=1}\frac{|\bar{x}_{i}-\bar{y}_{i}|^{4}}{4\epsilon^{4}}-\eta \bar{t}.
\end{equation}
Actually $\bar{x},\bar{y}$ and $\bar{t}$ depend on $\epsilon,\eta$, but we omit this dependence in the notation for  simplicity.\\
\textbf{2. Viscosity inequalities when  $\bar{t} > 0$.}  From the fundamental result of the User's guide to viscosity solutions \cite[Theorem 8.3]{GIL92}, for every $\rho >0$, we get  $a_{1},a_{2}\in \mathbb{R}$ and $\bar{X},\bar{Y}\in \mathscr{S}_{N}$ such that 
$$\dsps(a_{1},\bar{p},\bar{X})\in \bar{\mathcal{P}}^{2,+}u(\bar{x},\bar{t}),~~ \dsps(a_{2},\bar{p},\bar{Y})\in \bar{\mathcal{P}}^{2,-}v(\bar{y},\bar{t})$$
and 
$$
-\frac{1}{\rho}\left(
\begin{array}{cc}
I    &   0 \\

0    & I   \\
\end{array}
\right)
\leq \left(
\begin{array}{cc}
\bar{X}    &   0 \\

0    & -\bar{Y }  \\
\end{array}
\right)
\leq
\left(
\begin{array}{cc}
Z+\rho Z^{2}   &  -(Z+\rho Z^{2})   \\

-( Z+\rho Z^{2})   &   Z+\rho Z^{2}   \\
\end{array}
\right)
$$
$$ a_1-a_2= \eta,$$
for  $Z= D^{2}\varphi(\bar{x}-\bar{y})$, where $\varphi(x-y)= \displaystyle \sum^{N}_{i=0}\frac{|x_{i}-y_{i}|^{4}}{4\epsilon^{2}}$,\; $I$ is the identity matrix  and $$\dsps\bar{p}=\Big(\dsps\frac{|\bar{x}_{1}-\bar{y}_{1}|^{2}(\bar{x}_{1}-\bar{y}_{1})}{\epsilon^{4}},\cdot\cdot\cdot,\dsps\frac{|\bar{x}_{N}-\bar{y}_{N}|^{2}(\bar{x}_{N}-\bar{y}_{N})}{\epsilon^{4}}\Big).$$
Writing that $u$ is a subsolution and $v$ a supersolution of (\ref{HG}), we have 
\begin{equation}\label{vi}
\eta\leq\dsps  F_{\star}(\bar{x},\bar{t},\bar{p},\bar{X},K^{k,+}_{\bar{x},\bar{t},u})-F^{\star}(\bar{y},\bar{t},\bar{p},\bar{Y},{K}^{k,-}_{\bar{y},\bar{t},v}).\end{equation}

\begin{re}\rm{The reader should note a difficulty in applying \cite[Theorem 8.3]{GIL92} here. Indeed, one should double the time variable to prove \cite[Theorem 8.3]{GIL92}. It is not straightforward here because of the presence of the nonlocal term. This problem is solved by the stability result provided in \cite{SL03}. However, we give another proof of the stability in the appendix .}
\end{re}

\noindent\textbf{3. Comparison of the nonlocal terms .} From the definition of $M_{\epsilon,\eta}$,  for all $(x,y)$ $\in  \mathbb{R}^{N}\times \mathbb{R}^{N}$, we have 
\\
 $u(x,\bar{t})-v(y,\bar{t})-\displaystyle\sum^{N}_{i=1}\frac{|x_{i}-y_{i}|^{4}}{4\epsilon^{4}}\leq   u(\bar{x},\bar{t})-v(\bar{y},\bar{t})-\sum^{N}_{i=1}\frac{|\bar{x}_{i}-\bar{y}_{i}|^{4}}{4\epsilon^{4}}$.\\
Taking  $x=(\pi_{k}(\bar{x}),z)$ and  $y=(\pi_{k}(\bar{y}),z)$, for any $z$ $\in\mathbb{R}^{N-k}$, we have
$$u(\pi_{k}(\bar{x}),z,\bar{t})-v(\pi_{k}(\bar{y}),z,\bar{t})-\displaystyle\sum^{k}_{i=1}\frac{|\bar{x}_{i}-\bar{y}_{i}|^{4}}{4\epsilon^{4}} \leq u(\bar{x},\bar{t})-v(\bar{y},\bar{t})-\displaystyle\sum^{N}_{i=1}\frac{|\bar{x}_{i}-\bar{y}_{i}|^{4}}{4\epsilon^{4}}
,$$
thus we obtain 
\begin{equation}\label{comparaison}
u(\pi_{k}(\bar{x}),z,\bar{t})-u(\bar{x},\bar{t})\leq v(\pi_{k}(\bar{y}),z,\bar{t})-v(\bar{y},\bar{t})-\dsps\sum^{N}_{i=k+1}\frac{|\bar{x}_{i}-\bar{y}_{i}|^{4}}{4\epsilon^{4}}.
\end{equation}
\textbf{4. Upper-bound for the volume $K^{k,+}_{\bar{x},\bar{t},u}$and conclusion.}\\
 Since $u,v \in\mathcal{C}$ and by (\ref{maxp}) and (\ref{S}), $ \bar{x}$ and $\bar{y}$ remain  bounded independently of $\epsilon$ and $\eta$. Using $v$ $\in \mathcal{C}$, we  write  $$\dsps\mathop{\rm{lim}}_{\epsilon,\eta\rightarrow 0}{v(\bar{y},\bar{t})> -1},$$ 
then, there exist a positive constant $\mu$ independent  of $\epsilon$ and $\eta$, such that $v(\bar{y},\bar{t})\geq-1+\mu$. Then 
$$K^{k,+}_{\bar{y},\bar{t},v}\subset \{z\in \mathbb{R}^{N-k}: v(\pi_{k}(\bar{y}),z,\bar{t})\geq -1+\mu\}\subset \overline{B(0,R)}$$
where $R$ is a positive constant  independent of $\epsilon$ and $\eta$. By the inequality (\ref{comparaison}), if $ u(\pi_{k}(\bar{x}),z,\bar{t})-u(\bar{x},\bar{t})\geq 0$, we have  $v(\pi_{k}(\bar{y}),z,\bar{t})-v(\bar{y},\bar{t})\geq 0$, then    
\begin{equation}\label{ev}
K^{k,+}_{\bar{x},\bar{t},u}\subset K^{k,+}_{\bar{y},\bar{t},v} \subset\overline{B(0,R)} .
\end{equation}
 We distinguish  two  cases:\\
\textbf{First case}: $\dsps\bar{x}_{i}=\bar{y}_{i}$ for all $k < i\leq N$. In this case we have$$\dsps\bar{p}=\dsps\left(\frac{|\bar{x}_{1}-\dsps\bar{y}_{1}|^{2}(\bar{x}_{1}-\bar{y}_{1})}{\epsilon^{4}},\cdot\cdot\cdot,\dsps\frac{|\bar{x}_{k}-\bar{y}_{k}|^{2}(\bar{x}_{k}-\bar{y}_{k})}{\epsilon^{4}},0,\cdot\cdot\cdot,0\right).$$
In other words, $\tilde{\pi_{k}}(\bar{p})=0$. First, if $\bar{{x}}_{i}= \bar{{y }}_{i} $  for any $1\leq i \leq k$ then  $\bar{p}=0$. But, by \cite[Theorem 3.2]{GIL92}, and \cite[Lemma 2.4.3]{OBS}, there  exist  $\bar{X'},\bar{Y'}\in \mathscr{S}^{N}$ such that $$ \bar{X} \leq \bar{X'}\leq \bar{Y'}\leq \bar{Y},$$
with $\bar{X'}=\bar{Y'}=0$ when $\bar{p}=0$. Taking advantage of the ellipticity of $F^{\star}$ and $F_{\star}$, we get from (\ref{vi}) $$
\eta\leq\dsps  F_{\star}(\bar{x},\bar{t},0,{ O},K^{k,+}_{\bar{x},\bar{t},u})-F^{\star}(\bar{x},\bar{t},0,O,{K}^{k,-}_{\bar{y},\bar{t},v}).
$$
Assumption \textbf{(H2)} implies that  
\begin{center}
$\eta\leq 0,$
\end{center}
which is a contradiction since $\eta > 0$.\\
Therefore, there exists  $1\leq i_{0}\leq k$ such that   $\bar{{x}}_{i_{0}}\neq \bar{{y }}_{i_{0}}$ and $\bar{p}\neq 0$. Using \textbf{(H7)} and the viscosity inequality in $(\ref{vi})$, we get
$$
\eta\leq\dsps  G(\bar{x},\bar{t},\pi_{k}(\bar{p}),\bar{X})-G(\bar{y},\bar{t},\pi_{k}(\bar{p}),\bar{Y}).
$$
By \textbf{(H7)}, $G$ satisfies the assumptions \textbf{(H6-0)} then  
\begin{align}\label{ef}
\eta&\leq w\Big(|\bar{x}-\bar{y}|\left(1+|\pi_{k}(\bar{p})|\right)+\sigma|\bar{x}-\bar{y}|^{2}\Big)
\nonumber\\&\leq  w\Big(|\bar{x}-\bar{y}|\left(1+|\bar{p}|\right)+\sigma|\bar{x}-\bar{y}|^{2}\Big).
\end{align} 
Since $\dsps|\bar{p}|<\frac{|\bar{x}-\bar{y}|^{3}}{4\epsilon^{4}}$, then $\dsps|\bar{x}-\bar{y}||\bar{p}|\leq \frac{|\bar{x}-\bar{y}|^{4}}{4\epsilon^{4}}$. From this estimate, the inequality in (\ref{ef}) becomes 
\begin{equation}\label{iaf}
\eta\leq w\Big(|\bar{x}-\bar{y}|+\frac{|\bar{x}-\bar{y}|^{4}}{4\epsilon^{2}}+\sigma|\bar{x}-\bar{y}|^{2}\Big).
\end{equation}
Classical arguments \cite[Remark 3.8]{GIL92} show that $\dsps|\bar{x}-\bar{y}|$  and the penalisation term  $\dsps\sum^{N}_{i=1}\frac{|\bar{x}_{i}-\bar{y}_{i}|^{4}}{4\epsilon^{2}}$ tends to $0$ if $\epsilon$ tends to $0$, then $\dsps\frac{|\bar{x}-\bar{y}|^{4}}{4\epsilon^{2}}$ tends to $0$ if $\epsilon$ tends to zero. Finally we let  $\epsilon$  tend to 0 in $(\ref{ef})$, to obtain $\eta\leq 0$ which is a absurd.\\
\textbf{Second case}: There exists  $i_{0}$ such that $k< i_{0}\leq N$ and $\bar{x}_{i_{0}}\neq\bar{y}_{i_{0}}$. In this case, for all $z \in \mathbb{R}^{N-k}$, we have from $(\ref{comparaison})$ 
\begin{align*}
u(\pi_{k}(\bar{x}),z,\bar{t})-u(\bar{x},\bar{t})&\leq v(\pi_{k}(\bar{y}),z,\bar{t})-v(\bar{y},\bar{t})-\dsps\frac{|\bar{x}_{i_{0}}-\bar{y}_{i_{0}}|^{4}}{4\epsilon^{2}}
\\ &< v(\pi_{k}(\bar{y}),z,\bar{t})-v(\bar{y},\bar{t}).
\end{align*}
and therefore 
\begin{equation}\label{ctl}
K^{k,+}_{\bar {x},\bar{t},u} \subset{K}^{k,-}_{\bar{y},\bar{t},v}.
\end{equation}
Since $\bar{p}\neq 0$ in this case, the inequality viscosity in $(\ref{vi})$ becomes  
$$
\eta\leq\dsps  F(\bar{x},\bar{t},\bar{p},\bar{X},K^{k,+}_{\bar{x},\bar{t},u})-F(\bar{y},\bar{t},\bar{p},\bar{Y},{K}^{k,-}_{\bar{x},\bar{t},v}).
$$
From \textbf{(H3)} and (\ref{ctl}), we get 
$$
\eta\leq\dsps  F(\bar{x},\bar{t},\bar{p},\bar{X},K^{k,+}_{\bar{x},\bar{t},u})-F(\bar{y},\bar{t},\bar{p},\bar{Y},{K}^{k,+}_{\bar{x},\bar{t},u}).
$$
We use the assumption  \textbf{(H6-1)} to write  
\begin{equation}\label{es1}
\dsps \eta\leq w\bigg(|\bar{x}-\bar{y}|\Big(1+|\bar{p}|+\dsps\mathcal{L}^{N-k}(K^{k,+}_{\bar{x},\bar{t},u})^\frac{1}{N-k}|\tilde{\pi}_{k}(\bar{p})|+\sigma |\bar{x}-\bar{y}|^{2})|\Big)\bigg).
\end{equation}
Now, using  the estimate in  (\ref{ev}), we have 
$$\dsps\mathcal{L}^{N-k}(K^{k,+}_{\bar{x},\bar{t},u})^\frac{1}{N-k}\leq\mathcal{L}^{N-k}(\overline{B(0,R)})^{\frac{1}{N-k}}\leq C_{N-k}R$$
where  $C_{N-k}= \mathcal{L}^{N-k}\left(B_{N-k}\left(0,1\right)\right)^{\frac{1}{N-k}}$, and $R$ is independent of $\epsilon$. Then  
\begin{align}\label{estf}
\dsps\nonumber \eta&\leq w\bigg(|\bar{x}-\bar{y}|\Big(1+|\bar{p}|+ C_{N-k}R|\tilde{\pi}_{k}(\bar{p})|\Big)+\sigma|\bar{x}-\bar{y}|^{2}\bigg)
\\&\leq w\bigg(\dsps|\bar{x}-\bar{y}|+\frac{|\bar{x}-\bar{y}|^{4}}{4\epsilon^{2}}+ C_{N-k}R\dsps\frac{|\bar{x}-\bar{y}|^{4}}{4\epsilon^{2}}+\sigma|\bar{x}-\bar{y}|^{2}\bigg).
\end{align}
Since we know that $\dsps\frac{|\bar{x}-\bar{y}|^{4}}{4\epsilon^{2}}$ and  therefore $\dsps|\bar{x}-\bar{y}|$ converges to $0$ as $\epsilon$ goes to $0$, the above inequality (\ref{estf}) implies that  $\eta\leq 0$ which is a contradiction  since $ \eta > 0.$\\
\textbf{6. The case $\bar{t}=0.$} From the above step we obtain that the maximum $M_{\epsilon,\eta}$ is achieved for $\bar{t}=0$.  From (\ref{maxp}) and using the uniform continuity of $u_{0}$ ($u{_0}\in\mathcal{C}\cap C(\mathbb{R}^{N})$), for any $\rho \geq 0$, there exists $L_{\rho}>0$ such that
$$\frac{ M}{2}\leq M_{\epsilon,\eta} = u_{0}(\bar{x})-u_{0}(\bar{y})-\sum^{N}_{i=1}\frac{|\bar{x}_{i}-\bar{y}_{i}|^{4}}{4\epsilon^{2}}\leq \rho +L_{\rho}|\bar{x}-\bar{y}|\leq \rho +\frac{L_{\rho}^{2}\epsilon^{2}}{4},$$
which leads to a contradiction taking $\dsps\rho < \frac{M}{2}$ and sending $\epsilon$ to $ 0$. \hfill$\Box$


\section{Uniqueness result for non compact fronts}
\noindent In the previous section we considered the assumption \textbf{(H0)}, this assumption  forced us to deal with a compact initial front $\Gamma_{0}=\left\{u_{0}=0\right\}$. In this section we  deal with  non compact fronts. Instead of \textbf{(H0)}, we consider the assumptions: 

\vspace{3mm}

\noindent\textbf{(H0')} \textbf{Non compact initial  front }: $u_{0}\in BUC(\mathbb{R}^{N})$. 

\vspace{3mm}

\noindent{\bf(H1')} For all $(p,X)\in \mathbb{R}^{N}\times\mathscr{S}_{N}$, we have
$$\dsps\mathop{\rm{lim}}_{|p|,|X|\rightarrow 0}{F^{\star}(x,t,p,X,K^{k})}=\dsps\mathop{\rm{lim}}_{|p|,|X|\rightarrow 0}{F_{\star}(x,t,p,X,K^{k})}=0,$$
uniformly with respect to $(x,t,K^{k})\in \mathbb{R}^{N}\times[0,T]\times\mathscr{B}_{N-k}$.
\vspace{3mm}

\noindent\textbf{(H2')} For all $(x,t,K^{k})\in\mathbb{R}^{N}\times[0,T]\times\mathscr{B}_{N-k}$ and for all $(p,X,Y)$ in $\mathbb{R}^{N}\backslash\{0\}\times\mathscr{S}_{N}\times\mathscr{S}_{N}$ such that $(|p|+|X|),(|p|+|Y|)\leq R$, where $R$ is a positive constant, There exists a  nondecreasing modulus of continuity $w_{R}:\mathbb{R}_{+} \longrightarrow  \mathbb{R}_{+}$ and $w_{R}(0^{+})=0$  such that
$$F(x,t,p,X,K^{k})-F(x,t,q,Y,K^{k})\leq w_{R}\Big((|X-Y|_{\infty}+|p-q|)(1+|x|)\Big).$$  

\vspace{3mm}

\noindent\textbf{(H3')}  For any $K^{k}, L^{k}\in \mathscr{B}_{N-k}$ such that $K^{k}\backslash L^{k}\subset \mathbb{R}^{N-k}\backslash B(0,r)$, where $r>0$, and for bounded  $(p,X)\in\mathbb{R}^{N}\times\mathscr{S}_{N}$ we have
$$F(x,t,p,X,K^{k})-F(x,t,p,X,L^{k})\rightarrow 0$$
if $r$ tends to $+\infty$ uniformly with respect to $(x,t)\in\mathbb{R}^{N}\times[0,T]$ .

\vspace{3mm}

\noindent\textbf{(H4')} For any $K^{k}, L^{k}\in \mathscr{B}_{N-k}$, for any bounded  $(p,X)\in 
\mathbb{R}^{N}\backslash\{0\}\times\mathscr{S}_{N}$ such that $|\tilde{\pi_{k}}(p)|\leq \lambda$, we have 
$$F(x,t,p,X,K^{k})-F(x,t,p,X,L^{k})\rightarrow 0$$
if $\lambda$ tends to $0$ uniformly with respect to $(x,t)\in\mathbb{R}^{N}\times[0,T]$ .
\begin{re}
 \rm{The assumption \textbf{(H3') }implies that, if the Lebesgue measure of $K^{k}\backslash L^{k}$ is negligible, then the difference between $F(x,t,p,X,K^{k})$ and $ F(x,t,q,X,L^{k}) $ is small for any $x\in \mathbb{R}^{N}$ and  $X$ bounded and uniformly with respect to $t\in[0,T]$}. In addition, since we deal with the  non compact front \textbf{(H0')}, we  cannot control the nonlocal terms  and   instead of the assumption \textbf{(H5-1)}, we consider the assumption \textbf{(H5-0)}.
\end{re}
\begin{Th}\label{th2} Assume{\rm{ \textbf{(H0')- (H1')-(H2')-(H3')-(H4')-(H4)-(H5-0)-(H6-0)}}}. Let $u$ (resp. $v$) be a bounded upper-semicontinuous subsolution of (\ref{HG}) 
(resp.  bounded lower semicontinuous supersolution  of (\ref{HG})), then $u\leq$ $v$ in $\mathbb{R}^{N}\times [0,T]$.
\end{Th}
\begin{co}\label{CO2} Under the assumptions of Theorem $\ref{th2}$, there exists a unique viscosity solution $u$ of $(\ref{HG})$.
\end{co}
\noindent The  demonstration of this Corollary is postponed to (Section 5).

\vspace{3mm}

\noindent \textbf{Proof of the Theorem \ref{th2}}.\\
\textbf{1. The test-function.} We argue by contradiction. We suppose that there exists $(\tilde{x},\tilde{t})$ such that  
$$0<M=u(\tilde{x},\tilde{t})-v(\tilde{x},\tilde{t}).$$
In this case, we have to add some terms in the test-function in order to deal with non compact front. We consider 
\begin{equation}\label{max2}
\dsps M_{\epsilon,\eta,\alpha}=\displaystyle\mathop{\rm{sup}}_{\mathbb{R}^{N}\times\mathbb{R}^{N}\times\left[0,T\right]}\{u(x,t)-v(y,t)-(\sum^{N}_{i=1}\frac{\left|x_{i}-y_{i}\right|^{4}}{4\epsilon^{2}}+\alpha\left|x\right|^{2}+\alpha\left|
y\right|^{2})-\eta t\}.
\end{equation}
Since $u,v$ are bounded, for $\epsilon, \alpha,\eta >0$, the supremum is achieved at a point $(\bar{x},\bar{y},\bar{t})$   and  
for $\alpha, \eta$ are small enough, we have $$M_{\epsilon,\eta,\alpha}\geq \frac{M}{2}>0.$$
\noindent \textbf{2. Viscosity inequalities when $\bar{t}> 0$.} From the fundamental result of the User's guide to viscosity solutions \cite[Theorem 8.3]{GIL92}, for every $\rho >0$, we get  $a_{1},a_{2}\in \mathbb{R}$ and $\bar{X},\bar{Y}\in \mathscr{S}^{N}$ such that 
$$\dsps(a_{1},\bar{p}+2\alpha \bar{x},\bar{X}+2\alpha I)\in \bar{\mathcal{P}}^{2,+}u(\bar{x},\bar{t}),~~ \dsps(a_{2},\bar{p}-2\alpha \bar{y},\bar{Y}-2\alpha I)\in \bar{\mathcal{P}}^{2,-}v(\bar{y},\bar{t}),$$
$$ a_1 -a_2=\eta,$$

and 
$$
\frac{-1}{\rho}\left(
\begin{array}{cc}
I    &   0 \\

0    & I   \\
\end{array}
\right)
\leq \left(
\begin{array}{cc}
\bar{X}    &   0 \\

0    & -\bar{Y }  \\
\end{array}
\right)
\leq
\left(
\begin{array}{cc}
Z+\rho Z^{2}   &  -(Z+\rho Z^{2})   \\

-( Z+\rho Z^{2})   &   Z+\rho Z^{2}   \\
\end{array}
\right)
$$
for  $Z= D^{2}\varphi(\bar{x}-\bar{y})$, where $\varphi(x-y)= \displaystyle \sum^{N}_{i=0}\frac{|x_{i}-y_{i}|^{4}}{4\epsilon^{2}}$ for any $(x,y)\in\mathbb{R}^{N}\times\mathbb{R}^{N}$  and $$\dsps\bar{p}=\Big(\dsps\frac{|\bar{x}_{1}-\bar{y}_{1}|^{2}(\bar{x}_{1}-\bar{y}_{1})}{\epsilon^{4}},\cdot\cdot\cdot,\dsps\frac{|\bar{x}_{N}-\bar{y}_{N}|^{2}(\bar{x}_{N}-\bar{y}_{N})}{\epsilon^{4}}\Big).$$
Moreover, by \cite[Remark 3.8]{GIL92} and \cite[Proposition 2.5]{GGS}, $|\bar{X}|_{\infty}$,$|\bar{Y}|_{\infty}$ and  $|\bar{x}-\bar{y}|$ are bounded (independently of $\alpha$) and we have 
\begin{equation}\label{car}
\dsps\mathop{\rm{lim}}_{\epsilon\rightarrow 0} {\frac{|\bar{x}-\bar{y}|^{4}}{4\epsilon^{2}}}=0, \dsps\mathop{\rm{lim}}_{\epsilon,\alpha\rightarrow 0}{\alpha(|\bar{x}|^{2}+|\bar{y}|^{2})}=0,
\end{equation}
\begin{equation}\label{car1} 
\dsps\mathop{\rm{lim}}_{\alpha\rightarrow 0}{\alpha|\bar{y}|}=\dsps\mathop{\rm{lim}}_{\alpha\rightarrow 0}{\alpha|\bar{x}|}=0.
\end{equation}
Writing that $u$ is a subsolution and $v$ a supersolution of (\ref{HG}), we have 
\begin{equation}\label{vi2}
\eta\leq\dsps F_{\star}(\bar{x},\bar{t},\bar{p}+2\alpha\bar{x},\bar{X}+2\alpha I,K^{k,+}_{\bar{x},\bar{t},u})-F^{\star}(\bar{y},\bar{t},\bar{p}-\alpha\bar{y},\bar{Y}-2\alpha I,{K}^{k,-}_{\bar{y},\bar{t},v}).
\end{equation}
\textbf{3. Difference between $K^{k,+}_{\bar{x},\bar{t},u}$ and ${K}^{k,-}_{\bar{y},\bar{t},v}.$} 
 Because of the presence of  term $\alpha\left|x\right|^{2}+\alpha\left|
y\right|^{2}$  in the test function, the procedure used in the compact front case (Step 4) is non valid. Now, since $(\bar{x},\bar{y},\bar{t})$ is the maximum point of (\ref{max2}), we have
$$\begin{array}{l}
u(\pi_{k}(\bar{x}),z,\bar{t})-v(\pi_{k}(\bar{y}),z,\bar{t})\\

\hspace{2cm} 
\leq u(\bar{x},\bar{t})-v(\bar{y},\bar{t})-\dsps\sum^{N}_{i=k+1}\big(\frac{\left|\bar{x}_{i}-\bar{y}_{i}\right|^4}{4\epsilon^{4}}+\alpha (\bar{x}^{2}_{i}+\bar{y}^{2}_{i})\big)+2\alpha|z|^{2}.
\end{array}$$

From the definition of 
$K^{k,+}_{\bar{x},\bar{t},u}$ and $K^{k,-}_{\bar{y},\bar{t},v}$, we have $$K^{k,+}_{\bar{x},\bar{t},u}\subset K^{k,-}_{\bar{y},\bar{t},v}\cup\bar{\mathcal{E}},$$
where $\bar{\mathcal{E}}=K^{k,+}_{\bar{x},\bar{t},u}\cap\{v(\pi_{k}(\bar{y}),\cdot,\bar{t})\leq v(\bar{x},\bar{t})\}$. If $z\in \bar{\mathcal{E}}$, we have 
$$u(\bar{x},\bar{t})-v(\bar{y},\bar{t})\leq u(\pi_{k}(\bar{x}),z,\bar{t})-v(\pi_{k}(\bar{y}),z,\bar{t}).$$
It follows that 
\begin{equation}\label{E}
\bar{\mathcal{E}}\subset \left\{z\in \mathbb{R}^{N-k}:|z|^{2}\geq \frac{1}{2}\sum^{N}_{i=k+1}\Big(\frac{\left|\bar{x}_{i}-\bar{y}_{i}\right|^4}{4\alpha \epsilon^{4}}+ \bar{x}^{2}_{i}+\bar{y}^{2}_{i}\Big)\right\}.
\end{equation}
\textbf{4. Estimate of the term hand and side of the inequality $(\ref{vi2})$.} We distinguish two cases : 

\textbf{First case:} $\dsps\mathop{\rm{lim}}_{\alpha\rightarrow 0}{\bar{p}}=0$. By  \cite[Theorem 3.8]{GIL92} we have $\bar{X}\leq D^{2}_{\bar{x}\bar{x}}\varphi(\bar{x},\bar{y},\bar{t})$ and $\bar{Y}\geq -D^{2}_{\bar{y}\bar{y}}\varphi(\bar{x},\bar{y},\bar{t}).$
Taking advantage of the ellipticity of $F^{\star}$ and $F_{\star}$  in the viscosity inequality  (\ref{vi2}), we obtain
$$ 
\eta\leq\dsps F_{\star}(\bar{x},\bar{t},\bar{p}+2\alpha\bar{x},D^{2}_{\bar{x}\bar{x}}\varphi(\bar{x},\bar{y},\bar{t})+2\alpha I,K^{k,+}_{\bar{x},\bar{t},u})-F^{\star}(\bar{y},\bar{t},\bar{p}-\alpha\bar{y},-D^{2}_{\bar{y}\bar{y}}\varphi(\bar{x},\bar{y},\bar{t})-2\alpha I,{K}^{k,-}_{\bar{y},\bar{t},v}).$$
But, since $\dsps\mathop{\rm{lim}}_{\alpha\rightarrow 0}{\bar{p}}=0$, then $\dsps\mathop{\rm{lim}}_{\alpha\rightarrow 0}{D^{2}_{\bar{x}}\varphi(\bar{x},\bar{y},\bar{t})}=0$ and $\dsps\mathop{\rm{lim}}_{\alpha\rightarrow 0}{D^{2}_{\bar{y}}\varphi(\bar{x},\bar{y},\bar{t})}=0$. We send $\alpha$ to zero and we use  the assumption \textbf{(H1')} to write $\eta \leq 0$, which is a contradiction since $\eta  >0$.

\textbf{ Second case: $\dsps\mathop{\rm{lim}}_{\alpha\rightarrow 0}{\bar{p}}\neq 0$}. In which case, up to extract a subsequence there exists $\nu >0$ such that $\bar{p}\geq \nu$. Now, by (\ref{car1}) we can suppose that  $\bar{p}\neq 0, \bar{p}+2\alpha\bar{x}\neq 0$, $\bar{p}-2\alpha\bar{y}\neq 0$ and we have  the following estimate\\
$\noindent F(\bar{x},\bar{t},\bar{p}+2\alpha\bar{x},\bar{X}+2\alpha I,K^{k,+}_{\bar{x},\bar{t},u})- F(\bar{y},\bar{t},\bar{p}-2\alpha\bar{y},\bar{Y}-2\alpha I,K^{k,-}_{\bar{x},\bar{t},v})$\\
$\leq (I_{1})+(I_{2})+(I_{3})+(I_{4}),$\\
where\\
$(I_{1})= F(\bar{x},\bar{t},\bar{p}+2\alpha \bar{x},\bar{X}+2\alpha I,K^{k,+}_{\bar{x},\bar{t},u})-F(\bar{x},\bar{t},\bar{p},\bar{X},K^{k,+}_{\bar{x},\bar{t},u}),$\\
$(I_{2})=F(\bar{x},\bar{t},\bar{p},\bar{X} ,K^{k,+}_{\bar{x},\bar{t},u})-F(\bar{y},\bar{t},\bar{p},\bar{Y} ,K^{k,+}_{\bar{x},\bar{t},u}),$\\
$(I_{3})=F(\bar{y},\bar{t},\bar{p},\bar{Y} ,K^{k,+}_{\bar{x},\bar{t},u})-F(\bar{y},\bar{t},\bar{p},\bar{Y} ,K^{k,-}_{\bar{y},\bar{t},v}),$\\
$(I_{4})=F(\bar{y},\bar{t},\bar{p},\bar{Y} ,K^{k,-}_{\bar{y},\bar{t},v})-F(\bar{y},\bar{t},\bar{p}-2\alpha\bar{y},\bar{Y}-2\alpha I ,K^{k,-}_{\bar{y},\bar{t},v}).$\\
By the classical argument in \cite[Theorem 8.3]{GIL92} and \cite[Proposition 2.5]{GGS}, $\bar{p}+2\alpha\bar{x}, \bar{X}, \bar{Y}$ and $\bar{p}+2\alpha\bar{x}$ are bounded independently of $\alpha$. We suppose that $R_{\epsilon}={\rm{Max}}\Big((|\bar{p}+2\alpha\bar{x}|+|\bar{X}|),(|\bar{p}-2\alpha\bar{y}|+|\bar{Y}|)\Big)$, and we use the assumption 
\textbf{(H2')}, to obtain
\begin{align*}\label{es2.0}
(I_{1})&\leq w_{R_{\epsilon}}\Big((2\alpha| I|_{\infty}+2\alpha|\bar{x}|)(1+|\bar{x}|)\Big)\\&\leq w_{R_{\epsilon}} \Big(2\alpha+4\alpha|\bar{x}|+2\alpha|\bar{x}|^{2}\Big),
\end{align*}
since $|I|_{\infty}=1$. Moreover
\begin{align*}
(I_{4})&\leq w_{R_{\epsilon}}\Big((2\alpha| I|_{\infty}+2\alpha|\bar{y}|)(1+|\bar{y}|)\Big)\\&\leq w_{R_{\epsilon}} \Big(2\alpha+4\alpha|\bar{y}|+2\alpha|\bar{y}|^{2}\Big),
\end{align*}

Now, we use the assumption \textbf{(H6-0)}, to have
\begin{equation}\label{es2.2}
\begin{aligned}
(I_{2})&\leq w\Big(|\bar{x}-\bar{y}|(1+|\bar{p}|)+\sigma|\bar{x}-\bar{y}|^{2}\Big)
\\&\leq w\Big(|\bar{x}-\bar{y}|+ \frac{|\bar{x}-\bar{y}|^{4}}{4\epsilon^{2}}+\sigma|\bar{x}-\bar{y}|^{2}\Big).
\end{aligned}
\end{equation}
The estimate of the term $(I_{3})$ is prove later.\\
\noindent \textbf{5. End of the case  $\bar{t}>0$.}
By the above estimate, the viscosity inequality  in $(\ref{vi2})$ becomes
\begin{align}\label{es2.3}
\eta&\leq w_{R_{\epsilon}} \Big(2\alpha+4\alpha|\bar{x}|+2\alpha|\bar{x}|^{2}\Big)+w_{R_{\epsilon}} \Big(2\alpha+4\alpha|\bar{y}|+2\alpha|\bar{y}|^{2}\Big)
\\&+\nonumber w\Big(|\bar{x}-\bar{y}|+ \frac{|\bar{x}-\bar{y}|^{4}}{4\epsilon^{2}}+\sigma|\bar{x}-\bar{y}|^{2}\Big)+(I_{3}).
\end{align}
From (\ref{car}) and (\ref{car1}), we have 
\begin{equation}\label{es2.4}
\dsps\mathop{\rm{lim}}_{\epsilon,\alpha\rightarrow 0}{w_{R_{\epsilon}} \Big(2\alpha+4\alpha|\bar{x}|+2\alpha|\bar{x}|^{2}\Big)}=0.
\end{equation}
\begin{equation}\label{es2.7}
\dsps\mathop{\rm{lim}}_{\epsilon,\alpha\rightarrow 0}{w_{R_{\epsilon}}\Big(2\alpha+4\alpha^{2}|\bar{y}|+4\alpha^{2}|\bar{y}|^{2}\Big)}=0.
\end{equation}
\begin{equation}\label{es2.5}
\dsps\mathop{\rm{lim}}_{\epsilon,\alpha\rightarrow 0}{w\Big(|\bar{x}-\bar{y}|+ \frac{|\bar{x}-\bar{y}|^{4}}{4\epsilon^{2}}+\sigma|\bar{x}-\bar{y}|^{2}\Big)}=0.
\end{equation}
Now, we send $\alpha$ to zero, the inequality (\ref{es2.3}) becomes 
\begin{align}\label{es2.6}
\eta&\nonumber\leq \dsps\mathop{\rm{lim}}_{\alpha\rightarrow 0}\left({w_{R_{\epsilon}} \Big(2\alpha+4\alpha|\bar{x}|+2\alpha|\bar{x}|^{2}\Big)}+ {w_{R_{\epsilon}} \Big(2\alpha|I|+4\alpha|\bar{y}|+2\alpha|\bar{y}|^{2}\Big)}\right)
\\&+\dsps\mathop{\rm{lim}}_{\alpha\rightarrow 0}{ w\Big(|\bar{x}-\bar{y}|+ \frac{|\bar{x}-\bar{y}|^{4}}{4\epsilon^{2}}+\sigma|\bar{x}-\bar{y}|^{2}\Big)}+
 \dsps\mathop{\rm{lim}}_{\alpha\rightarrow 0}(I_{3}).\end{align}
From now on, we denote by 
$\dsps I_{\epsilon,\alpha}$ the sum$\frac{1}{2}\sum^{N}_{i=k+1}\frac{\left|\bar{x}_{i}-\bar{y}_{i}\right|^4}{4\epsilon^{4}}$. Then (\ref{E}) becomes
\begin{equation}
\bar{\mathcal{E}}\subset \left\{z\in \mathbb{R}^{N-k}:|z|^{2}\geq \dsps I_{\epsilon,\alpha}+\sum^{N}_{i=k+1} (\bar{x}^{2}_{i}+\bar{y}^{2}_{i})\right\}.
\end{equation}
For $(I_3)$ we distinguish two cases :\\
\textbf{First case}: $\dsps \mathop{\rm{lim}}_{\alpha \rightarrow 0}{I_{\epsilon,\alpha}}=0$, in this case we have $\dsps\mathop{\rm{lim}}_{\alpha\rightarrow 0}{|\tilde{\pi}_{k}(\bar{p})|}=0$. Since $\bar{Y}-2\alpha I$ and $\bar{p}-2\alpha\bar{y}$ remain  bounded (independently of $\alpha$), we use  the assumption \textbf{(H4')}, to obtain
$$\dsps\mathop{\rm{lim}}_{\alpha\rightarrow 0}{(I_3)}=0.$$
 Now, we send $\epsilon$ to zero in (\ref{es2.6}) and we use (\ref{es2.4}),(\ref{es2.7}) and (\ref{es2.5}), the inequality (\ref{es2.6})
becomes $\eta\leq 0$, which is a contradiction since $\eta >0$.\\
\textbf{Second case:} $\dsps \mathop{\rm{lim}}_{\alpha \rightarrow 0}{I_{\epsilon,\alpha}}\neq 0$, in this case, up to extract a subsequence there exists $\delta>0$ such that   $$\dsps I_{\epsilon,\alpha}\geq \delta.$$ 
The estimate in  (\ref{E}) implies   
 $$\bar{\mathcal{E}}=K^{k,+}_{\bar{x},\bar{t},u}\backslash K^{k,-}_{\bar{y},\bar{t},v}\subset\mathbb{R}^{N-k}\backslash B(0,r_{\alpha}).$$
where  $\dsps r_{\alpha}=\dsps\frac{\dsps I_{\epsilon,\alpha}}{\alpha}\geq \frac{\delta}{\alpha}$ tend to $+\infty$  if $\alpha$ tends to zero. Since we know that  $\bar{Y}-2\alpha I$ and $\bar{p}-2\alpha\bar{y}$ are  bounded (independently of $\alpha$), then the assumption \textbf{(H3')} implies
$\mathop{\rm{lim}}_{\alpha\rightarrow 0}{(I_3)}=0$, and we obtain a contradiction as in the first case.\\
\textbf{End of the proof.}  From above, we have necessarily $\bar{t}=0$, and we conclude as in the precedent section, Step $6$. \hfill$\Box$
\section{Existence result} 
 In this Section, we use the classical  Perron's method to prove the proof of Corollary \ref{CO}, when  the initial front is  compact since the one in the non compact case (Corollary \ref{CO2}) can be adapted easily. This proof is given by  the  three following steps.\\
\textbf{Step 1.} In  this step we construct a subsolution $\underline{u} \in\mathcal{C}$ and a supersolution $\overline{u}\in \mathcal{C}$ of our equation  $(\ref{HG})$. We start with the following Lemma.
\begin{Le}\label{L0} Let $ A\leq -\Big( 3 L_{1}+ 2L_{N-k}L_{1}+2(L_{2}+1)(2N+3)\Big)$, where $L_{1}, L_{2}$ appear in the assumption (H5) and  $C_{N-k}= \mathcal{L}^{N-k}\left(B\left(0,1\right)\right)^{\frac{1}{N-k}}$. Then the function $$g(x,t)=\frac{e^{At}}{1+|x|^{2}}-1$$   is a smooth subsolution of (\ref{HG})  for $t>0$.
\end{Le}
The proof of the Lemma \ref{L0} is postponed. Now, from the subsolution $g$ for $t>0$, we construct a subsolution   in the class $\mathcal{C}$ which satisfy the initial condition. We consider the nondecreasing  function $\phi$ from $(-\infty,0)$ to $\mathbb{R}$ defined by $$\dsps\phi(t)=\mathop{{\rm{inf}}}_{y\in \mathbb{R}^{N}}\{u_{0}(y): g(y,0)\geq t\}.$$
By the definition of $\phi$ we have \begin{equation}\label{ci} 
\phi(g(y,0))\leq u_{0}(y).
\end{equation}
Starting from this function $\phi$, we  will build a regular function which has the property (\ref{ci}). For this reason, we need to prove  the following Lemmas.
\begin{Le}\label{L2} The function $\phi_{n}:(-\infty,\frac{1}{n})\rightarrow \mathbb{R}$ defined by $$ \dsps\phi_{n}(t)=n \int^{t}_{t-\frac{1}{n}}\phi(s)ds.$$  have the following properties:
\\ 
i)- $\phi_{n}(t)\leq \phi(t)$ for all $t\in (-\infty,0)$.\\
ii)-   For all $n$ the function $t \rightarrow\phi_{n}(t)$ is nondecreasing.\\
iii)- For all $n$ the function $\phi_{n}\in C((0,\frac{1}{n}))$ and  $\phi_{n}(g) \in \mathcal{C}$.
\end{Le} 
\noindent \textbf{Proof of the Lemma 4.2.}\\
i)- Since $\phi$ is nondecreasing we have $\phi(s)\leq\phi(t)$ for all $s \leq t$ and     
$$ \dsps\phi_{n}(t)=n \int^{t}_{t-\frac{1}{n}}\phi(s)ds \leq n \int^{t}_{t-\frac{1}{n}}\phi(t)ds\leq \phi(t)$$ 
then $\phi_{n}(t)\leq \phi(t)$ for all $t\in (-\infty,0).$\\
ii)- Let $t^{'}\geq t$, by  change of the variables $s^{'}=s+(t^{'}-t)$, we have 
\begin{align*}
\dsps\phi_{n}(t)= n\int^{t}_{t-\frac{1}{n}}\phi(s)ds=&\; n \int^{t^{'}}_{t^{'}-\frac{1}{n}}\phi\Big(s^{'}-(t^{'}-t)\Big)ds^{'}\\&\leq n \int^{t^{'}}_{t^{'}-\frac{1}{n}}\phi(s')ds^{'}\\&=\phi_{n}(t^{'})
\end{align*}
thus $\phi_{n}$ is nondecreasing.\\
iii)- The function $\phi$ is bounded, then the function $\phi_{n}$ is $C((-\infty,\frac{1}{n}))$. Since $\phi_{n}$ is nondecreasing by  ii) we have 
$$\mathop{\rm{lim}}_{|x|\rightarrow\infty}{\phi_{n}\Big(g(x,t)\Big)}=\mathop{\rm{inf}}_{x\in\mathbb{R}^{N}}{\phi_{n} \Big( g(x,t)\Big)} 
=\phi_{n}(-1)=-1$$ 
and $\phi_{n} \Big( g(x,t)\Big)> -1$ for any $(x,t)\in \mathbb{R}^{N}\times[0,T]$. Then $\phi_{n}(g) \in \mathcal{C}.$\hfill$\Box$\\
Now we consider the function $${\tilde{\phi}}_{n}(t)=n \int^{t}_{t-\frac{1}{n}}\phi_{n}(t)ds.$$
Since $\phi_{n}\in C((-\infty,\frac{1}{n}))$,  then $\tilde{\phi}_{n} \in C^{1}((-\infty,\frac{1}{n}))$ and the Lemma $\ref{L2}$ remains valid with the function $\tilde{\phi}_{n}$ (we replace $\phi$ by $\phi_{n}$). We argue in the same manner we consider the function $$\dsps\hat{\phi}_{n}(t)=n \int^{t}_{t-\frac{1}{n}}\tilde{\phi}_{n}(t)ds.$$
Since $\tilde{\phi}_{n}$ is continuous, then $\hat{\phi}_{n}\in C^{1}(-\infty,\frac{1}{n})$ and satisfies the following Lemma.
\begin{Le}\label{L1} The function $\hat{\phi}_{n}$ have the following properties \\
i)- $\hat{\phi}_{n}(t)\leq \phi(t)$ for all $t\in (-\infty,0)$.\\
ii)-  For all $n$ the function $t \rightarrow\hat{\phi}_{n}(t)$ is nondecreasing.\\
iii)-  For all $n$ the function $\hat{\phi}_{n}\in C^{2}((-\infty,\frac{1}{n}))$ and  $\hat{\phi}_{n}(g) \in \mathcal{C}$.
\end{Le}
\noindent The proof of this Lemma is similar to the proof of Lemma \ref{L2}.\\
We use the assumption \textbf{(H4)} which ensures that, the fronts is invariant by nondecreasing changes $u \rightarrow \psi(u)$(see \cite{ES91,GGC91}). Since  $g$ is a subsolution of (\ref{HG}) for $t>0$ (Lemma $\ref{L0})$, and   $\hat{\phi}_{n}$ is nondecreasing, the function $\hat{{\phi}}_{n}\circ g$
is a subsolution of $(\ref{HG})$ for $t>0$. The definition of $\phi$ and the Lemma \ref{L2} implies
$$  \hat{\phi}_{n}\Big(g(x,0)\Big)\leq\phi(g(x,0))\leq u_{0}(x)~~  {\rm{and}}~~ \hat{\phi}_{n}(g) \in \mathcal{C},$$
then $\hat{\phi}_{n}\circ g$ is a subsolution of (\ref{HG}) for all $t\in [0,T]$.

For the construction of the supersolution, we argue in the same manner but, we start  in the Lemma $\ref{L0}$ with  $\dsps f(x,t)=\frac{e^{A_{0}t}}{1+|x|^{2}}-1$, where $A_{0}\geq \Big(3L_{1}+2C_{N-k}L_{1}+2(2N+3)(L_{2}+1)\Big)$, and instead of  the function $\phi$, we take the  nondecreasing function 
$$ \dsps\psi(t)=\mathop{{\rm{sup}}}_{y\in \mathbb{R}^{N}}\{u_{0}(y): f(y,0)\leq t\}.$$
Then, there exists a suite of   nondecreasing functions $\hat{\psi}_{n}\in C^{1}(-\infty,\frac{1}{n})$ such that   $\hat{\psi}_{n}\circ f \in \mathcal{C}$ is a supersolution of $(\ref{HG})$. \hfill$\Box$   

Now, it is enough to take $\underline{u}=\hat{\phi}_{n}(g)$ and $\overline{u}=\hat{\psi}_{n}(f)$  to conclude  the proof of  Step 1 .\\
\textbf{Step 2.} Consider the set $\mathcal{F}$ of subsolution of (\ref{HG}) $w$ such that $\underline{u}\leq w\leq \overline{u}$. Set then for every $(x,t)\in \mathbb{R}^{N}\times[0,T], v(x,t)=\dsps\mathop{\rm{sup}}_{w\in\mathcal{F}}{w(x,t)}$. By Step 1, the set $\mathcal{F}$ is nonempty and $v$ is well-defined. Thus, we get from the comparison result and classical arguments of the Perron's method that $v$ is a discontinuous solution of (\ref{HG}). For the proof of these classical arguments we refer the reader to Crandall, Ishii, Lions \cite{GIL92} and to  Barles \cite{B94}.\\
Now, as  the subsolutions and the sursolutions do not satisfy the condition initial with equality, 
 we need to using  the some arguments used in \cite[Proposition 1]{AT96}  and \cite [Theorem 4.7]{B94} to conclude that $v^{\star}(\cdot,0)=v_{\star}(\cdot,0)=u_{0}$. Thus we deduce from the comparison result that $v^{\star}=v_{\star}=v$ which is the desired continuous solution.\\
\textbf{Step 3.} We show that the solution $v$ built in Step 2 is actually  in $\mathcal{C}$. If $w\in \mathcal{F}$, then $w\in \mathcal{C}$ and $w\leq \bar{u}$. Since  $\bar{u}\in \mathcal{C}$ we conclude easily   $v\in \mathcal{C}$.\hfill$\Box$\\
\noindent\textbf{Proof of the Lemma \ref{L0}.} First, an easy computation show
$$\dsps|D^{2}g|_{\infty} \leq \frac{2e^{At}(2N+3)}{(1+|x|^{2})^{2}}.$$
By \textbf{(H5)}, to show that $g$ is  subsolution, it suffices to prove  be  that
\begin{eqnarray*}
(E)&=&\dsps\frac{\partial g}{\partial t}(x,t)+L_{1}\Big((1+|x|)|Dg(x,t)|\Big)|\tilde{\pi}_{k}(Dg(x,t))|\\
& &+L_{1}\Big(\mathcal{L}^{N-k}(K^{k,+}_{x,t,g})^{\frac{1}{N-k}}|\tilde{\pi}_{k}(Dg(x,t))|\Big)
+(L_{2}+1)|x|^{2}|D^{2}g(x,t)| 
\end{eqnarray*}
is non positive.\\
Since  $g$  is radial and decreasing in $|x|$, we have  $$K^{k,+}_{x,t,g}\subset B(0,|x|)~~{\rm{and}}~ \mathcal{L}^{N-k}(K^{k,+}_{x,t,g})\leq C_{N-k}|x|^{N-k}$$
or equivalently  
\begin{equation}\label{v}
\mathcal{L}^{N-k}(K^{k,+}_{x,t,g})^{\frac{1}{N-k}}\leq C_{N-k}|x|.\end{equation}
Using $(\ref{v})$, $|\tilde{\pi}_{k}(Dg(x,t))|\leq |Dg(x,t)|$, and the above estimate of the upper bound of $|D^{2}g|_{\infty}$, we have 
\begin{equation}\label{A}(E)\leq \frac{Ae^{At}}{1+|x|^{2}}+L_{1}\Big((1+|x|)\frac{2e^{At}|x|}{(1+|x|^{2})^{2}}+\frac{2e^{At }C_{N-k}|x|^{2}}{(1+|x|^{2})^{2}}\Big)+(L_{2}+1)\frac{2e^{At}(2N+3)}{(1+|x|^{2})^{2}}|x|^{2}.
\end{equation}
We develop the right-hand side term  of $(\ref{A})$, we obtain
$$(E)\leq \frac{e^{At}(A+3L_{1}+2C_{N-k}L_{1}+2(L_{2}+1)(2N+3))}{1+|x|^{2}}.$$
Since $A\leq -\Big(3L_{1}+2C_{N-k}L_{1}+2(2N+3)(L_{2}+1)\Big)$, then $(E)\leq 0$, and $g$ is a subsolution of $(\ref{HG})$ for $t>0$ . $\hspace{7.2cm}\Box$ 

\section{Appendix}
In this Section, we prove a stability result  for our nonlocal equation. In fact, we need this result in the  proof of  Theorem \ref{th1}, Theorem \ref{th2} and the Perron's method (See Remark 2.1). In \cite {SL03}, this result was formulated in the following way. 
\begin{Th}\label{sr} If $(u_{n})_{n\geq 1}$ is a sequence  of upper-semicontinuous viscosity subsolutions  of  
\begin{equation}\label{subFn}
\displaystyle\frac{\partial u}{\partial t}(x,t)+ F_{n}(x,t,Du,D^{2}u,K^{k,+}_{{x,t,u}})=0,
\end{equation}
(resp. lower-semicontinuous supersolutions) of  
\begin{equation}\label{supFn}
\displaystyle\frac{\partial u}{\partial t}(x,t)+ F_{n}(x,t,Du,D^{2}u,K^{k,-}_{{x,t,u}})=0,
\end{equation}
where $(F_{n})_{n\geq1}$ is a sequence of uniformly local bounded functions on $\mathbb{R}^{N}\times[0,T]\times\mathbb{R}^{N}\backslash\{0\}\times\mathscr{S}_{N}\times\mathscr{B}_{N-k}$ with satisfying  the monotonicity condition \textbf{{\rm{(H3)}}}. We suppose that, the functions $(u_{n})_{n\geq 1}$ are uniformly local bounded on $\mathbb{R}^{N}\times[0,T]$, then $\bar{u}=\dsps\mathop{{\rm{\lim sup^{\star}}}}_{n}{u_{n}}$ {\rm{(}}resp. $\underline{u}=\dsps\mathop{{\rm{\lim inf_{\star}}}}_{n}u_{n})$ is a subsolution (resp. supersolution) of
$$\frac{\partial u}{\partial t}+\underline{F}(x,t,Du,D^{2}u,K^{k,+}_{{x,t,u}})=0,$$ 
respectively of
$$\frac{\partial u}{\partial t}+\overline{F}(x,t,Du,D^{2}u,K^{k,-}_{{x,t,u}})=0,$$
where 
$$\underline{F}(x,t,p,X,K^{k})=\dsps\mathop{{\rm{\lim inf_{\star}}}}_{n}F_{n}(x_{n},t_{n},p_{n},X_{n},K^{k}_{n}),$$ 
respectively
 $$\underline{F}(x,t,p,X,K^{k})=\dsps\mathop{\rm{\lim sup ^{\star}}}_{n}F_{n}(x_{n},t_{n},p_{n},X_{n},K^{k}_{n}),$$ 
when $ x_{n}\rightarrow $  $x$, $p_{n}\rightarrow p$,$X_{n}\rightarrow X$,$t_{n}\rightarrow t$, and $ K^{k}_{n}$ converge to $K^{k}$ if $n$ tends to $\infty$ .
\end{Th}
\noindent\textbf{Proof of the Theorem \ref{sr}}: 
We give the proof only for $\overline{u}$, that for $\underline{u}$ being similar. Let $(x_{0},t_{0})\in\mathbb{R}^{N}\times[0,T]$ and  $\phi\in C^{\infty}(\mathbb{R}^{N}\times[0,T])$ such that $\overline{u}-\phi$ has maximum at $(x_{0},t_{0})$. Then from the definition of $\overline{u}$ follows that there is a subsequence of $(u_{n})_{n\geq1}$ that we also denote by $(u_{n})_{n\geq1}$ such that $u_{n}-\phi$ has  maximum at $(x_{n},t_{n})$ and  
$$(x_{n},t_{n})\rightarrow (x_{0},t_{0})\quad {\rm{and}}\quad u_{n}(x_{n},t_{n})\rightarrow \overline{u}(x_{0},t_{0})\quad \rm{as} \quad n\rightarrow \infty.$$\\
So since $u_{n}$ are subsolutions of (\ref{subFn}), we have 
\begin{equation}\label{sriv}
\frac{\partial \phi}{\partial t}+F_{n}(x_{n},t_{n},D\phi(x_{n},t_{n}),D^{2}\phi(x_{n},t_{n}),K^{k,+}_{{x_{n},t_{n},u_{n}}})\leq0.
\end{equation}
To continue the proof we need the following Lemma.
\begin{Le}\label{ls} Let $u=\dsps\mathop{{{\rm{\lim sup^{\star}}}}}_{n}{u_{n}}$ where  $(u_{n})_{n\geq1}$ is a sequence of uniformly local bounded  and upper semicontinuous functions from $\mathbb{R}^{N}\times[0,T]$ to $\mathbb{R}$  and $(x_{n},t_{n})$  a sequence of\;   $\mathbb{R}^{N}\times[0,T]$ such that  $(x_{n},t_{n})\rightarrow (x,t)$ and $u_{n}(x_{n},t_{n}) \rightarrow u(x,t)$  as   $n\rightarrow \infty$, then 
$$\dsps\mathop{\rm{\lim sup}}_{n}{\1_{K^{k,+}_{x_{n},t_{n},u_{n}}}(z)}\leq \1_{K^{k,+}_{{x,t,u}}}(z),$$
for all $z\in\mathbb{R}^{N}$. Or equivalently  we have 
$$\dsps\mathop{\limsup}_{n}{K^{k,+}_{x_{n},t_{n},u_{n}}}\subset K^{k,+}_{x,t,u},$$
for the set convergence in the topology defined in the introduction. 
\end{Le}
\noindent The proof of the Lemma \ref{ls} is postponed. Now, we return to the proof of Theorem \ref{sr}. We use the Lemma \ref{sr}  to  conclude that 

$$\dsps\mathop{\rm{\lim sup}}_{n}{K^{k,+}_{x_{n},t_{n},u_{n}}\cup{K^{k,+}_{x,t,\bar{u}}}}=K^{k,+}_{x,t,\bar{u}},$$
then   $\1_{K^{k,+}_{x_{n},t_{n},u_{n}}\cup K^{k,+}_{{x_{0},t_{0},\bar{u}}}}$ converges  to $\1_{K^{k,+}_{{x_{0},t_{0},\bar{u}}}}$ in L$^{1}_{loc}(\mathbb{R}^{N-k})$ as $n\rightarrow\infty$.  By the monotonicity condition of $F_{n}$, we use (\ref{sriv}) to write  
$$\frac{\partial \phi}{\partial t}+F_{n}(x_{n},t_{n},D\phi(x_{n},t_{n}),D^{2}\phi(x_{n},t_{n}),K^{k,+}_{{x_{n},t_{n},u_{n}}}\cup K^{k,+}_{{x_{0},t_{0},\bar{u}}})\leq0.$$
Moreover, the regularity  of the test-function $\phi$ enables us to affirm that
$$
\dsps\frac{\partial \phi}{\partial t}(x_{0},t_{0})+\underline{F}(x_{0},t_{0},D\phi(x_{0},t_{0}),D^{2}\phi(x_{0},t_{0}), K^{k,+}_{{x_{0},t_{0},\bar{u}}})$$
$$\hspace{2cm}\leq \frac{\partial \phi}{\partial t}(x_{0},t_{0})+\dsps\mathop{{\rm{\lim inf}}}_{n }{F_{n}(x_{n},t_{n},D\phi(x_{n},t_{n}),D^{2}\phi(x_{n},t_{n}),K^{k,+}_{{x_{n},t_{n},u_{n}}}}\cup K^{k,+}_{{x_{0},t_{0},\bar{u}}})$$
$\hspace{2cm}\leq 0.$\\
Then, $\bar{u}$ is a subsolution of(\ref{subFn}). \hfill$\Box$ 

\vspace{3mm}

\noindent\textbf{Proof of Lemma \ref{ls}.} For the proof of the Lemma \ref{ls} we distinguish two cases:\\
\textbf{First case}: $u(z,t)\geq u(x,t)$, in this case  we have $\1_{{K^{k,+}_{{x,t,u}}}}(z)=1$ and the result is clear  since, obviously  $$\dsps\mathop{\rm{\lim sup}}_{n}{\1_{K^{k,+}_{x_{n},t_{n},u_{n}}}(z)}\leq \1_{K^{k,+}_{{x,t,u}}}(z)=1.$$
\textbf{Second case:} $u(z,t)<u(x,t)$, in this case since $u_{n}(x_{n},t)\rightarrow u(x,t) $  and by the definition of $u $ for $n$ large enough we have $u_{n}(z,t_{n})<u_{n}(x_{n},t_{n})$ and $$\dsps\mathop{\limsup}_{n}{\1_{K^{k,+}_{x_{n},t_{n},u_{n}}}(z)}\leq \1_{K^{k,+}_{{x,t,u}}}(z).$$ \hfill$\Box$
\begin{re}\rm{ We point out that Lemma \ref{ls} is not true if we replace ``$>$'' by  ``$\geq$'' in the definition of $K^{k,+}$. It explain why we need to change ``test-sets''for supersolutions in Definition \ref{ds} (see Remark \ref{rs})}. \hfill$\Box$
\end{re}  

\noindent\textbf{Acknowledgments.}\\
The author would like  to thank Guy Barles  and Olivier Ley  for  their great help, their support, their  enriching discussions  and their many fruitful suggestions in the preparation of this article.  I also would like to express my great fullness to  Romain Abraham and  Maïtine Bergounioux for introducing me to the tomography reconstruction and the studied model. This work was supported by grants of the center region and the national center of scientific research CNRS.

\end{document}